\documentclass[12pt]{amsart}
\usepackage{amsmath,amssymb}

\newtheorem{thm}{Theorem}
\newtheorem*{thm2}{Theorem}
\newtheorem{lem}[thm]{Lemma}

\theoremstyle{definition}

\theoremstyle{remark}

\theoremstyle{plain}

\newcommand{\Q}{{\mathbf Q}}

\newcommand{\C}{{\mathbf C}}

\newcommand{\G}{{G(f,\psi)}}
\newcommand{\Gt}{{G(f,\psi_t)}}

\newcommand{\zp}{\zeta_p}
\newcommand{\zn}{\zeta_{n}}
\newcommand{\zpk}{\zeta_{p^k}}
\newcommand{\zpl}{\zeta_{p^l}}

\newcommand{\PK}{{{\mathfrak P}_K}}
\newcommand{\PL}{{{\mathfrak P}_L}}
\newcommand{\PP}{{{\mathfrak P}}}
\newcommand{\OL}{{{\mathfrak O}_L}}
\newcommand{\OK}{{{\mathfrak O}_K}}
\newcommand{\OO}{{{\mathfrak O}}}

\newcommand{\gal}{\operatorname{Gal}}

\newcommand{\pexp}{\eta}

\numberwithin{equation}{section}

\title{On a character sum problem of H. Cohn}
\author{P\"ar Kurlberg}
\address{Department of Mathematics, University of Georgia, Athens GA
  30602 ({\tt kurlberg@math.uga.edu})}

\thanks{Author supported in part by the National Science Foundation
  (DMS 0071503).}

\begin{document}

\hyphenation{gran-ville}

\begin{abstract}
  Let $f$ be a complex valued function on a finite field $F$ such that
  $f(0) = 0$, $f(1) = 1$, and $|f(x)| = 1$ for $x \neq 0$. Cohn asked
  if it follows that $f$ is a nontrivial multiplicative character
  provided that $\sum_{x \in F} f(x) \overline{f(x+h)} = -1$ for
  $h \neq 0$.
  We prove that this is the case for finite fields of prime
  cardinality under the assumption that the nonzero values of $f$ are
  roots of unity.
\end{abstract}

\maketitle

\section{Introduction}

Let $p$ be prime and let $F_{p^k}$ be the finite field with $p^k$
elements. Let $f : F_{p^k}^\times \rightarrow \C$ be a nontrivial
multiplicative character, and extend $f$ to a function on $F_{p^k}$ by
letting $f(0)=0$. It is then easy to see that the following holds:
\begin{equation}
\label{e:cohn}
\sum_{x \in F_{p^k}}
f(x) \overline{f(x+h)} =
\begin{cases} 
-1 & \text{if $h \neq 0$}\\ p^k-1 & \text{if $h = 0$}
\end{cases}
\end{equation}

Cohn asked (see p. 202 in \cite{montgomery-book}) if the converse is
true in the following sense: if a function $f : F_{p^k} \rightarrow
\C$  satisfies 
\begin{equation}
\label{e:cohn-extra}
f(0) = 0, \ f(1) = 1, \ \text{and } |f(x)| = 1 \text{ for $x\neq0$}
\end{equation}
and equation \ref{e:cohn}, does it follow that $f$ is a multiplicative
character?

The problem has recently received some attention. In
\cite{false-cohn}, Choi and Siu proved that the converse is not true
for $k>1$.  One of the arguments given is quite pretty, and proceeds
as follows: Let $\lambda$ be a linear automorphism of $F_{p^k}$ so
that $\lambda(1)=1$. If $f$ satisfies \ref{e:cohn} and
\ref{e:cohn-extra}, so does $f$ composed with $\lambda$.  Now, if $f$
is an injective multiplicative character then the converse being true
implies that $f$ composed with $\lambda$ must be an injective
multiplicative character.  On the other hand, a simple counting
argument shows that the number of possible $\lambda$'s is greater than
the number of injective characters.

However, the case $k=1$ remains unresolved. In \cite{biro-cohn}, Biro
proved that there are only finitely many functions satisfying equation
\ref{e:cohn} and \ref{e:cohn-extra} for each $p$. Biro also solved the
following ``characteristic $p$'' version of the problem
(\cite{biro-cohn}, Theorem 2):

\begin{thm2}[Biro]
\label{t:biro}
Let $p$ be a prime, let $F_p$ be the finite field with $p$ elements,
and 
$F \supset F_p$ any field of characteristic $p$. Assume that there is
given an $a_i \in F$ for every $i \in F_p$ such that $a_0=0, a_1=1,
a_i \neq 0$ for $i \neq 0$, and 
$$
\sum_{i \in F_p^\times} \frac{a_{i+j}}{a_i} =  -1
$$
for every $j \in F_p^\times$. Then $a_i=i^A$ for every $i \in F_p$
with some $1\leq A \leq p-2$. 
\end{thm2}

Using this Biro deduces that the converse holds for functions taking
values in $\{-1, 0, 1\}$.\footnote{There appears to be several
  independent proofs of this result, see the introduction in
  \cite{false-cohn}.}  In fact, if $m$ is coprime to $p$, then the
case of the nonzero values of $f$ being $m$-th roots of unity can be
deduced in a similar way: Let $\OO$ be the ring of integers in
$\Q(e^{2 \pi i/m})$, and let $\PP \subset \OO$ be a prime ideal lying
above $p$. The result then follows from the theorem by letting $F =
\OO/\PP$ and noting that $m$-th roots of unity are distinct modulo $p$.
(Since $|f(x)|=1$ for $x \neq 0$ we have $\overline{f(x)}=1/f(x)$.)

The aim of this paper is to show that the converse is true for the
case $k=1$,
under the 
additional assumption that the nonzero values of $f : F_p \to \C$ are
$m$-th roots of unity, including the case $p|m$.
We begin by giving a proof that does not depend on Biro's result for
the case $(m,p)=1$, and we then show how to modify the argument for
the general case.

{\em Acknowledgements:} I would like to thank Ernest Croot, Andrew
Granville, Robert Rumely, and Mark Watkins for helpful and stimulating
discussions.  
I would also like to thank the referee for several suggestions on how
to improve the exposition, and for pointing out that the case $p|m$
can be deduced independently of Biro's theorem.



\section{Preliminaries}

In what follows we assume that $p$ is odd since the case $p=2$ is
trivial. 

We will use the following conventions: if a function $f$ takes values
in $\C$ and $\sigma \in \operatorname{Aut}(\C/\Q)$, then we let
$f^\sigma$ be the function defined by $f^\sigma(x)=\sigma(f(x))$. We
regard $\psi(x)=e^{2\pi i x/p}$ as a nontrivial additive character of
$F_p$. For an integer $t$, $\psi_t$ will denote the character
$\psi_t(x)=\psi(tx)$. By $\zeta_m$ we denote the $m$-th root of unity
$\zeta_m = e^{2 \pi i/m}$. 

Let $m$ be even and large enough so that all nonzero values of $f$ are
$m$-th roots of unity, and write $m=np^k$, where $(n,p)=1$. Let
$K=\Q(\zeta_n)$, $L=K(\zp,\zpk)$, and let $G=\gal(L/\Q)$,
$H=\gal(L/K)$ denote the Galois groups of $L/\Q$ and $L/K$. By $\OK$
and $\OL$ we will denote the ring of integers in $K$ respectively $L$.



The ``Gauss sum'' 
$$
\G = \sum_{x=0}^{p-1} f(x) \psi(x)
$$
is clearly an algebraic integer. As in the case of classical
Gauss sums, the absolute value of $\G$ can easily be determined:

\begin{lem}
\label{l:G-abs-value}
  If $f$ satisfies \ref{e:cohn}, then
$$
|\Gt|=
\begin{cases}
\sqrt{p} & \text{if $t \not \equiv 0 \mod p$,} \\
0 & \text{if $t \equiv 0 \mod p$.}
\end{cases}
$$  

\end{lem}
\begin{proof} We have
$$
|\Gt|^2 = \sum_{x,y \in F_p} f(x) \overline{f(y)} \psi(t(x-y)) =
\sum_{x,h \in F_p} f(x) \overline{f(x+h)} \psi(-th) 
$$
$$
= \psi(0) \sum_{x \in F_p} f(x) \overline{f(x)} +
\sum_{h \in F_p^\times} \psi(-th) 
\sum_{x \in F_p} f(x) \overline{f(x+h)} 
$$
$$
=
p-1 
- \sum_{h \in F_p^\times} \psi(-th) 
= 
\begin{cases}
p & \text{if $t \not \equiv 0 \mod p$,} \\
0 & \text{if $t \equiv 0 \mod p$,} 
\end{cases}
$$

\end{proof}

The action of complex conjugation on $K$ is given by an element in
$G$, and since $G$ is abelian, equation \ref{e:cohn} is $G$-invariant.
I.e., if $f$ satisfies \ref{e:cohn-extra}, so does $f^\sigma$ for all
$\sigma \in G$. But if $\sigma \in G$ then
$\sigma(G(f,\psi))=G(f^{\sigma},\psi_t)$, where
$\sigma(\zeta_p)=\zeta_p^t$. Since $f^{\sigma}$ also satisfies
\ref{e:cohn}, we find that $|G(f^{\sigma},\psi_t)|=p^{1/2}$, and hence
the $\Q$-norm of $\G$ is a power of $p$. The factorization of the
principal ideal $\G\OL$ thus consists only of prime ideals $\PL|p$.

It is well known 
that $\Q(\zpk)/\Q$ is totally ramified over $p$, and that
$\Q(\zeta_n)/\Q$ does not ramify at $p$ if $(n,p)=1$.  Comparing
ramification indices gives that if $\PK$ is a prime ideal in $\OK$
that divides $p$, then $\PK$ is totally ramified in $L$. In
particular, if $\PL$ is any prime ideal in the ring of integers in
$\OL$ that lies above $p$, then $\sigma(\PL)=\PL$ for all $\sigma \in
H$.

Let $l = \max(1,k)$. Then $H$ consists of elements $\sigma_t$ such
that  
\begin{eqnarray*}
\sigma_t(\zpl)=\zpl^t, &  
\sigma_t(\zeta_n)=\zeta_n.
\end{eqnarray*}
Choose $t$ so that $\sigma_t$ generates $H$.  
Applying $\sigma_t$ to the
principal ideal 
$$\G \OL =\prod_{\PL|p} \PL^{\pexp(\PL)}$$
we find that
$$\sigma_t(\G \OL) = 
\sigma_t( \prod_{\PL|p} \PL^{\pexp(\PL)} ) =
\prod_{\PL|p} \PL^{\pexp(\PL)} =
\G \OL 
$$
and hence $\sigma_t(\G)=u\G$ for some unit $u$. 

Since the absolute value of any complex embedding of $\G$ equals
$\sqrt{p}$, we find that all conjugates of $u=\sigma(\G)/\G$ has
absolute value one. Hence $u$ is in fact a root of unity, and there
are integers $a,b$ such that

\begin{equation}
\label{e:g-mult}
\sigma_t(\G)  = \zpl^a \zeta_n^b \G.
\end{equation}

\section{The case $(m,p)=1$}

Since $f$ is fixed by $H$ we find that $\sigma_t(\G)=\Gt$, and 
equation~\ref{e:g-mult} can, after the change of variable $x \to
t^{-1}x$,  be written as
\begin{equation}
\label{e:sum1}
\sum_{x=1}^{p-1} f(x) \psi(x)
=
\zeta_p^{-a} \zeta_n^{-b}
\sum_{x=1}^{p-1} f(t^{-1}x) \psi(x).
\end{equation}

\begin{lem}
  If $f$ takes values in $n$-th roots of unity for $x \not \equiv 0 \mod
  p$ and equation \ref{e:sum1} holds then $a \equiv 0 \mod p$.
\end{lem}
\begin{proof}
From \ref{e:sum1} we obtain that 
\begin{equation}
\label{e:sum2}
\sum_{i=1}^{p-1} A_i \zeta_p^i
=
\sum_{i=0}^{p-1} B_i \zeta_p^{i}
\end{equation}
where $A_i=f(i)$ and $B_i=\zeta_n^{-b} f(t^{-1}(i+a))$. (Note that
$B_{p-a}=0$.) Since $ 1= -\sum_{i=1}^{p-1} \zeta_p^i$
we may rewrite
\ref{e:sum2} as 
\begin{equation}
\sum_{i=1}^{p-1} A_i \zeta_p^i
=
\sum_{i=1}^{p-1} (B_i-B_0) \zeta_p^{i}.
\end{equation}
The elements $\{\zeta_p, \zeta_p^2, \zeta_p^3, \ldots
\zeta_p^{p-1}\}$ are
linearly independent over $K$, hence
$A_i=B_i-B_0$. From lemma~\ref{l:G-abs-value} we have
$\sum_{x=0}^{p-1} f(x)=0$, which implies that
$\sum_{i=1}^{p-1} A_i = 0$, as well as $\sum_{i=0}^{p-1} B_i = 0$.
Therefore,
$$
0 =
\sum_{i=1}^{p-1} A_i =
\sum_{i=1}^{p-1} (B_i-B_0) =
\sum_{i=0}^{p-1} B_i-pB_0 = -pB_0.
$$
But $B_0=\zeta_n^{-b} f(t^{-1}(0+a))$ which is nonzero unless $a \equiv 0
\mod p$. 
\end{proof}

Thus
\begin{equation}
\label{e:sum3}
\sum_{x=1}^{p-1} f(x) \psi(x)
=
\zeta_n^{-b}
\sum_{x=1}^{p-1} f(t^{-1}x) \psi(x)
\end{equation}
and the linear independence of $\{\zeta_p, \zeta_p^2,
\zeta_p^3, \ldots \zeta_p^{p-1}\}$ over $K$ implies that
$$
f(t^{-1}x)=f(x) \zeta_n^{b }
$$
for all $x \neq 0$. Thus
$$
f(t^{-k})= 
f(t^{-(k-1)}) \zeta_n^{b} = \ldots =
f(1) \zeta_n^{kb} =\zeta_n^{kb}.
$$
Taking $k=p-1$ we find that $\zn^{b}$ is a $(p-1)$-th root of unity, 
and that $f$ is a multiplicative character.

\section{The general case}
In this case $m=np^k$ where $(n,p)=1$ and $k>0$. We will need the
following:

\begin{lem} 
\label{l:zp-coeff}
If $a_i \in K$ and  $\sum_{i=0}^{p^k-1} a_i \zpk^i \in K(\zp)$ then 
\begin{equation}
\label{e:in-K-zp}
\sum_{i=0}^{p^k-1} a_i \zpk^i =
\sum_{j=0}^{p-1} a_{p^{k-1}j} \zp^j 
\end{equation}
\end{lem}
\begin{proof}
We may assume that $k>1$. The minimal polynomial for $\zpk$ (over $K$
as well as over $\Q$) is given by 
$$
\frac{x^{p^k}-1}{x^{p^{k-1}}-1} = 1 + x^{p^{k-1}} + x^{2p^{k-1}}
+ \ldots + x^{(p-1)p^{k-1}}.
$$
Hence, by letting $\tilde{i} \in [0,p^{k-1}-1]$ be a representative
of $i$ modulo $p^{k-1}$, we can rewrite the left hand side of equation
\ref{e:in-K-zp} as
$$
\sum_{i=0}^{(p-1)p^{k-1}-1} (a_i-a_{(p-1)p^{k-1}+\tilde{i}}) \zpk^i 
$$
with no further relations among the $\zpk^{i}$'s, and thus
$$
\sum_{i=0}^{(p-1)p^{k-1}-1} (a_i-a_{(p-1)p^{k-1} +\tilde{i}}) \zpk^i
\in K(\zp)
$$
if and only if $a_i-a_{(p-1)p^{k-1}+\tilde{i}}=0$ for all $i$ not
congruent to zero modulo $p^{k-1}$.
\end{proof}


Recall from equation~\ref{e:g-mult} (note that $l=k$ since $k\geq 1$)
that  
$$
\sigma_t(\G)= \zpk^a \zn^b \G.
$$
Let $\tilde{G}=\zpk^s \G$ where
$\sigma_t(\zpk^s)/\zpk^s = \zpk^{-a}$. (Such an $s$ exists as
$\sigma_t(\zpk^s)/\zpk^s = \zpk^{(t-1)s}$, and $t \not \equiv 1
\mod p$ since $\sigma_t$ generates $H$.) We then have
$$
\sigma_t(\tilde{G}) = 
\sigma_t(\zpk^s \G) 
$$
$$
= \sigma_t(\zpk^s) \sigma_t( \G) =
\sigma_t(\zpk^s) \zpk^a \zn^b \G) =
\zn^b \tilde{G}.
$$ 

The following lemma shows that $\tilde{G}$ must transform by a nontrivial
$n$-th root of unity:

\begin{lem} 
\label{l:nontrivial-n-root}
  There is no integer $s$ such that $\zpk^s \G \in K$.
\end{lem}

\begin{proof}
  We first assume that $\zpk^s=1$.  Let $\G \OL = \prod_{\PL|p}
  \PL^{\pexp(\PL)}$ be the factorization of the principal ideal $\G
  \OL$.  Since $p$ does not ramify in $K$, we have $p \OK =
  \prod_{\PK|p} \PK$, and hence $p \OL = \prod_{\PL|p} \PL^e$  where
  $e$ is the ramification index of $\PK$ in $L$.

Since $\psi(x) = \zp^x$ is congruent to $1$ modulo $\PL$ for all $x$,
we find that 
$$
\G=
\sum_{x=0}^{p-1} f(x) \psi(x) 
\equiv 
\sum_{x=1}^{p-1} f(x) \mod \PL
$$
for all $\PL|p$. Now, since $f(0)=0$, we have $\sum_{x=1}^{p-1}
f(x) =  G(f,\psi_0)$ and by
lemma~\ref{l:G-abs-value}, $G(f,\psi_0) = 0$.  Thus $\G \in \PL$ for
all $\PL|p$, i.e., $\pexp(\PL)>0$ for all $\PL|p$.
But if $\G \in K$ then $e|\pexp(\PL)$ for all $\PL|p$, and since
complex conjugation permutes the set of primes of $\OL$ that lies
above $p$, and
$$
p=\G \overline{\G},
$$
we get that $\PL^{2e}| p \OL $ for all $\PL$, contradicting that the
ramification index is $e$.

For the general case, the previous argument carries through by noting
that $\zpk^s$ is a unit (and thus multiplication of $\G$ by $\zpk^s$
does not change the ideal factorization) and that $\G \in \PL$ if and
only if $\zpk \G \in \PL$.
\end{proof}

Since $\sigma_t$ has order $p^{k-1}(p-1)$ and $(n,p)=1$ we find that
$\zn^b$ must  be a nontrivial $(p-1)$-th root of unity. Hence
there exists a nontrivial multiplicative character $\chi$ of
$F_p^\times$ such that $\chi(t^{-1})=\zeta_n^b$. But $\sigma_t(
G(\chi,\psi))=G(\chi,\psi_t)= \chi(t^{-1}) G(\chi,\psi)$ and thus
$$
\delta =\frac{\tilde{G}}{G(\chi,\psi)}
$$
is $\sigma_t$-invariant and hence an element of $K$. Moreover,
$|\delta|=1$ (for all complex embeddings) since
$|\tilde{G}|=|G(\chi,\psi)|=p^{1/2}$. 

Write $f(x)=f_1(x)f_2(x)$ where $f_1(x)$ takes values
in $p^k$-th roots of unity and $f_2(x)$ takes values in $n$-th roots of
unity. We will show that $f_1(x)$ must be constant. 

\begin{lem}
Let 
$$
a_i=\sum_{x: \zpk^s f_1(x)\psi(x) = \zpk^i} f_2(x)
$$
If 
\begin{equation}
\label{e:f-delta-chi}
\zpk^s \sum_{x=1}^{p-1} f(x) \psi(x) =
\delta \sum_{x=1}^{p-1} \chi(x) \psi(x), 
\end{equation}
then $|a_{i}|=0$ unless $i=p^{k-1}j$ for $j=1,2,\ldots, p-1$, in which case 
$|a_{i}|=1$. 
In particular, $\zpk^s f_1(x) \psi(x)$ ranges over all nontrivial
$p$-th roots of unity. 

\end{lem}
\begin{proof}

Collecting terms in \ref{e:f-delta-chi} according to the values of
$\zpk^s f_1(x)\psi(x)$, we obtain
\begin{equation}
\label{e:zeta-rel1}
\sum_{i=0}^{p^k-1} a_i \zpk^i
=
\delta \sum_{i=1}^{p-1} \chi(i) \zp^i \in K(\zp).
\end{equation}

Clearly $a_i \in K$ and $a_i \neq 0$ for at most $p-1$ values of
$i$. 
%
Letting $A_i=a_{p^{k-1}i}$ we may, by lemma \ref{l:zp-coeff}, write
equation \ref{e:zeta-rel1} as 
$$
\sum_{i=0}^{p-1} A_i \zp^i
=
\delta \sum_{i=1}^{p-1} \chi(i) \zp^i.
$$
Since $1= -\sum_{i=1}^{p-1} \zeta_p^i$ we get that
$$
\sum_{i=1}^{p-1} (A_i-A_0) \zp^i
=
\sum_{i=0}^{p-1} A_i \zp^i
=
\delta \sum_{i=1}^{p-1} \chi(i) \zp^i
$$
and hence $A_i-A_0=\delta \chi(i)$ for all $i$. 

Since $a_i \neq 0$ for at most $p-1$ values of $i$, $A_0 \neq 0$
implies that $A_j=0$ 
for some $j \neq 0$, and thus $|A_0|=|\delta \chi(j)-A_j|=1$. Since
$$
0=
\delta \sum_{i=1}^{p-1} \chi(i) 
=
\sum_{i=1}^{p-1} (A_i-A_0) 
=
\sum_{i=0}^{p-1} A_i -  pA_0,
$$
we find that $|\sum_{i=0}^{p-1} A_i| = p|A_0|=p$. On the other
hand, $|\sum_{i=0}^{p-1} A_i| \leq \sum_{x=1}^{p-1} |f_2(x)|=p-1$.
Thus $A_0=0$, and it follows that $A_i=\delta \chi(i)$ for $i \neq 0$.
In other words, $a_{p^{k-1}j} = A_j =\delta \chi(j)$ for $j = 1, 2,
\ldots, p-1$, and since there are at most $p-1$ nonzero values among
the $a_i$'s, the remaining ones must all be equal to zero.
\end{proof}

Now, the lemma gives that $\zpk^s f_1(1) \psi(1) = \zpk^s \zp$ is a
$p$-th root of unity, hence $p^{k-1}$ must divide $s$, and the nonzero
values of $f_1(x) \psi(x)$ are thus distinct $p$-th roots of unity.
Replacing 
$\psi$ by $\psi_r$, for $r \not \equiv 0 \mod p$, in the previous
argument gives that $f_1(x) \psi(rx)$ also ranges over distinct
$p$-th roots of unity. On the other hand, if $f_1(x)$ is not constant,
then there exists $r \not \equiv 0 \mod p$ such that the set $\{
f_1(x) \psi_r(x) \}_{x=1}^{p-1}$ contains strictly less than $p-1$
elements. (If $f_1(x_1) \neq f_1(x_2)$, write $f_1(x_1)=\zp^{y_1},
f_1(x_2)=\zp^{y_2}$ and take $r \equiv -(y_2-y_1)(x_2-x_1)^{-1} \mod
p$.)  Hence $f_1(x)$ must be constant, and since $f_1(1)=1$, we find
that the nonzero values of $f(x)$ are in fact $n$-th roots of unity.
The result has thus been reduced to the case $(m,p)=1$.

\bibliographystyle{amsplain}

\end{document}